%% file: FeinbergSAFinalA.tex
  \documentclass[11pt]{article}

\usepackage[utf8]{inputenc}

\usepackage[square,numbers]{natbib}
\usepackage{fullpage}
\usepackage{xcolor}
\usepackage{amsmath}
\usepackage{amssymb}
\usepackage{bold-extra}
\usepackage{mathtools}
\usepackage{setspace}
\usepackage{enumitem}
\usepackage{amsthm}
\newcommand{\slim} {\mathop{\rm lim\,sup\,}}
\newcommand{\ilim} {\mathop{\rm lim\,inf\,}}

\def\U{\mathbb{U}}

\def\S{\mathbb{S}}
\def\K{\mathbb{K}}

\def\Y{\mathbb{Y}}
\def\Yb{\mathbb{Y}}

\def\X{\mathbb{X}}
\def\Xb{\mathbb{X}}
\newcommand{\Xc}{\mathcal{X}}
\newcommand{\Hc}{\mathcal{H}}

\def\A{\mathbb{A}}
\def\Ab{\mathbb{A}}
\def\Rb{\mathbb{R}}

\def\R{\mathbb{R}}
\def\P{\mathbb{P}}
\def\F{\mathbb{F}}
\def\T{\mathcal{T}}

\def\B{\mathcal{B}}

\DeclareMathOperator{\diag}{\mathbf{Diag}}
\newtheorem{lemma}{Lemma}[section]
\newtheorem{theorem}[lemma]{Theorem}
\newtheorem{corollary}[lemma]{Corollary}
\newtheorem{definition}[lemma]{Definition}

\begin{document}

\title{On Pioneering Works of Albert Shiryaev on Markov Decision Processes and Some Later Developments}
\author{Eugene A. Feinberg\thanks{%
        Department of Applied Mathematics and Statistics,  Stony Brook University, Stony Brook, NY 11794-3600, USA, (e-mail: {eugene.feinberg@stonybrook.edu}).
   } }
%
%
\maketitle

\abstract{This article is dedicated to three fundamental papers on Markov Decision Processes and on control with incomplete observations published by Albert Shiryaev approximately sixty years ago.  One of these papers was coauthored with O.V. Viskov.  We discuss some of the results and some of many rich ideas presented in these papers and survey some later developments.  At the end we mention some recent studies of Albert Shiryaev on Kolmogorov's equations for jump Markov processes and on control of continuous-time jump Markov processes.}

\section{Introduction}
\label{sec:1}
Albert Shiryaev is one of the pioneers and major creators of the theory of controlled stochastic processes.  His contributions are deep and broad.  They include sequential analysis, statistics of stochastic processes, change point problems, theory of martingales, limit theorems for stochastic processes, stochastic differential equations, mathematical finance, and many other fields.   This article describes some results and research directions in the theory of Markov Decision Processes (MDPs) influenced by and related to Shiryaev's papers \cite{Sh1,Sh2,VS} on MDPs and on control of stochastic processes with incomplete information.

MDPs deal with control of stochastic processes.  This field provides mathematical foundations to Reinforcement Learning \cite{BT,SB}, which is one of the major areas of Artificial Intelligence.  For example, introduced by the DeepMind team in 2017 the famous program AlphaZero~\cite{Silver}  played chess better than other computer programs and humans.  In addition, it was learning how to play chess within approximately 9 hours of self-training. The program demonstrated the similar success with several other difficult games.

Viskov and Shiryaev~\cite{VS} wrote a foundational paper on infinite-horizon MDPs with average costs per unit time.  Shiryaev~\cite{Sh1,Sh2} studied  problem with incomplete information in discrete and continuous time.  The corresponding now popular model in MDPs is a Partially Observable Markov Decision Process (POMDP), which is one of the major models in Reinforcement Learning.  The analysis of POMDPs is based on their reduction to MDPs with states being probability distributions of states of the original POMDP.  This reduction was formulated in \cite{Sh1,Sh2}.

There are two important questions in the theory of MDPs: (i) what is the structure of optimal and nearly optimal policies, and (ii) how to compute such policies.  The central facts are the validity of optimality equations sometimes called Bellman equations, and the possibility to reduce  MDP problems to linear programming problems.

\section{MDPs with Finite State and Action Sets}
\label{sec:2}
An MDP is defined by the set of states $\X,$ set of actions $\A,$ transition probability $p,$ and one-step cost function $c.$  In this section $\X$ and $\A$ are finite sets. The time parameter $t=0,1,\ldots$ is discrete.  If at a state $x\in\X$ an action $a\in\A$ is selected, then the process moves to the next state $z\in\X$ with probability $p(z|x,a),$ and the one-step cost $c(x,a)$ is incurred. Costs can be either real-valued or equal to $+\infty.$

In most of the studies, action sets at different states can be different.  That is, $A(x)\subset\A$ is the set of actions available at the state $x\in\X.$ Here we do not consider the sets $A(x)$ because we allow the possibility  $c(x,a)=+\infty$ for some $x\in\X$ and $a\in\A.$  If action sets $A(x)$ are state-dependent, where $x\in\X,$ it is possible to set $c(x,a)=+\infty$ for $a\in \A\setminus A(x)$ and make action sets state-independent.

Let $H:=(\X\times\A)^\infty$ be the set of all trajectories, and $H_t:=X\times (\A\times\X),$, $t=0,1,\ldots,$ be the sets of histories.  A policy $\pi$ is defined as a sequence $(\pi_0,\pi_1,\ldots)$ of conditional probabilities on $\A$ given $h_t\in H_t,$ where $h_t=x_0,a_0,x_1,a_1,\ldots,x_t$ is the history by epoch  $t=0,1,\ldots,$ where $x_t$ and $a_t$ are the state and actions at the epoch $t=0,1,\ldots.$  At time $t$ the next action $a_t$ is selected by the distribution $\pi_t(\cdot|h_t).$  If a policy $\pi$ always chooses actions with probabilities 0 or 1, it is called nonrandomized.  If the choice of an action depends only on state and time, the policy is called randomized Markov.  Such a policy is called Markov if it is nonrandomized.  If the choice of an action depends only on the current state, the policy is called stationary.  A nonrandomized stationary policy is called deterministic. A deterministic policy is defined by a function $\phi:\X\to\A.$ Let $\Pi$ be the set of all policies and $\F$ be the set of all deterministic policies.

A policy $\pi$ and an initial state $x$ define a probability on the set of trajectories $H_\infty,$ and this probability is denoted by $P_x^\pi.$  This standard fact follows from the Ionescu Tulcea theorem and also from the Kolmogorov extension theorem. Expectations with respect to this probability are denoted by $E_x^\pi.$

For infinite-horizon problems, the standard two objective criteria are the expected total discounted costs
\begin{equation}\label{eqdicd}
v_\beta^\pi(x)=E_x^\pi \sum_{t=0}^\infty \alpha^t c(x_t,a_t),
\end{equation}
where $\alpha\in [0,1)$ is the discount factor,  and average costs per unit time
\begin{equation}
w^\pi(x)=\limsup_{t\to\infty} \frac{1}{T}E_x^\pi \sum_{t=0}^{T-1} c(x_t,a_t).
\end{equation}

In general, if the objective function is $g^\pi(x),$ then the value function is $g(x):=\inf_{\pi\in\Pi} g^\pi(x),$ where $x\in\X.$  A policy $\pi$ is called optimal if $g^\pi(x)=g(x)$ for all $x\in\X.$  For $\epsilon>0,$ a policy $\pi$ is called $\epsilon$-optimal if $g^\pi(x)\le g(x)+\epsilon$ for all $x\in\X.$

Shapley~\cite{Shapley} introduced stochastic games with expected total discounted costs and proved the existence of equilibrium stationary policies for zero-sum stochastic games with expected total discounted payoffs.  This is essentially a more general fact than the existence of deterministic optimal policies for discounted MDPs with finite state and action sets. However, the direct proof for MDPs \cite{Blac,VS} is easier, and it follows from the Banach fixed point theorem.

Blackwell~\cite{Blac}, Derman~\cite{Der}, and Viskov and Shiryaev~\cite{VS} independently proved the existence of deterministic optimal policies for average-cost MDPs.  Blackwell~\cite{Blac} did not formulate this fact.  He proved that there exist $\alpha^*\in [0,1)$ and a deterministic policy $\phi$ such that $\phi$ is discount-optimal for all discount factors $\alpha\in [\alpha^*,1).$ Such policies are called Blackwell-optimal now. Blackwell-optimal policies are average-cost optimal, but an average-cost optimal policy may not be Blackwell-optimal; see Puterman~\cite[Example 10.1.1] {Pu}.

The proofs in \cite{Der,VS} are based on the Tauberian theorem: for a sequence $(b_t)_{t=0}^\infty,$ let us consider the Cezaro and Abel lower and upper limits
\begin{equation}
C_*:=\liminf_{T\to\infty}\frac{1}{T} \sum_{t=0}^{T-1}b_t\quad {\rm and}\quad  C^*:=\limsup_{T\to\infty}\frac{1}{T} \sum_{t=0}^{T-1}b_t,
\end{equation}
\begin{equation}
A_*:=\liminf_{\alpha\uparrow\infty}(1-\alpha)\sum_{t=0}^\infty \alpha^tb_t\quad {\rm and}\quad A^*:= \limsup_{\alpha\uparrow\infty}(1-\alpha)\sum_{t=0}^\infty \alpha^tb_t;
\end{equation}
then
\begin{equation}
C_*\le A_*\le A^*\le C^*.
\end{equation}
The relevant beautiful fact is the Hardy-Littlewood theorem stating that $A_*= A^*$ implies $C_*= C^*.$ An example in \cite{BFZ} shows that it is possible that $C_*= A_*$ and $A^* = C^*,$ but $A_*< A^*.$

Optimality equations play an important role in the theory of MDPs.  For $x\in\X$ and $a\in\A,$ we define $P^af(x):=\sum_{z\in\X}  p(z|x,a)f(z),$ where $f:\X\to\mathbb{R}\cup{-\infty}. $ In other words, $P^a f(x):=E\{ f(x_1)|x_0=x, a_0=a\}, $  and in this form this definition holds for problems with infinite state and action sets under minimal measurability and integrability assumptions.  Let us also define $T_\beta^af(x):=c(x,a)+\beta P^a f(x),$ where $\beta\ge 0.$ Then the optimality operator is \[
T_\beta f(x):=\min_{a\in\A} T_\beta^a f(x).
\]
This definition holds in general if minimum is replaced with infimum. We also define operators $T_\beta^\phi f(x)\mapsto T_\beta^{\phi(x)} f(x),$ where $\phi\in\F$ is a deterministic policy.  We usually write $T$ instead of $T_\beta$ if $\beta=1.$

Then  $v_\beta^\phi=T_\beta^\phi v^\phi$ for a problem with the discount factor $\beta\in [0,1),$  the optimality equation
\begin{equation}\label{eqoptD}
v_\beta=T_\beta v_\beta
\end{equation}
holds, and these two equations have the unique solutions $v_\beta^\phi$ and $v_\beta$ respectively.  A deterministic policy is optimal for a discounted MDP if and only if for all $x\in\X$
\begin{equation}
\phi(x)\in A^*(x):=\{a\in\A:\, v_\beta(x)= T_\beta^a v_\beta(x)\}.
\end{equation}

For average-cost MDPs with finite state and action sets, for each deterministic policy $\phi,$ the system of equations

\begin{equation}\label{eqCan}
\left\{
\begin{aligned}
w^\phi &=& P^\phi w^\phi  \\
w^\phi+u^\phi &=& T^\phi u^\phi
\end{aligned}
\right. ,
\end{equation}
with two unknown  variables $w^\phi$ and $u^\phi$ uniquely defines $w^\phi.$  For $\phi\in\F$ and for $u^\phi$ satisfying \eqref{eqCan},  if $w^\phi=Pw^\phi$ and
$w^\phi+u^\phi=Tu^\phi,$ then the policy is called canonical. A canonical policy exists for an MDP with finite sets of states and actions, and a canonical policy is average-cost optimal~\cite{DY}.
%
%

Value and policy iteration algorithms are two main methods for solving discounted MDPs. We do not formulate them here since they are broadly known and used; see, e.g., Puterman~\cite{Pu}.   The value iteration algorithm is based on iterating the right-hand size of optimality equation~\eqref{eqoptD}. Also,    optimality equation~\eqref{eqoptD} can be used to write a linear program (LP), and the policy iteration algorithm implements the simplex method with the block-pivoting rule applied to the dual LP; see, e.g.,  \cite{DF,Kal,Kal1,Pu}.  There is also a version of the policy iteration algorithm implementing the simplex method with Dantzig's pivoting rule, but usually it is slower.

Because of the link to LPs, finding an optimal policy for and MDP is a weakly polynomial problem.  In addition, Tseng~\cite{Ts} proved that value iterations are weakly polynomial. Ye~\cite{Ye} discovered that policy iterations are strongly polynomial if the discount factor is fixed and viewed as a constant.  This was a remarkable discovery in linear programming extended to other LPs in Kitahara and Mizuno~\cite{KM}.  For discounted MDPs Scherrer~\cite{Scher} improved some of Ye's~\cite{Ye} estimates.  Post and Ye~\cite{PY} proved that the policy iteration algorithm with Dantzig's pivoting rule is strongly polynomial for deterministic MDPs for all discount factors $\beta\in [0,1).$

An example in \cite{FHu} demonstrates that that value iterations are not strongly polynomial.  However, value iterations guarantee exponentially fast convergence of value functions, and value iterations are strongly polynomial for computing $\epsilon$-optimal policies \cite{FHe}.

For average-cost MDPs, policy iterations are usually used.  Value iterations are also possible under some assumptions; Federgruen and Schweitzer~\cite{FedS}.

Also, if $w^\phi(x)$ is  constant in $x$ for all $\phi\in \F,$ then the first equation in \eqref{eqCan} always holds, and the second equation leads to the optimality  equation
\begin{equation}\label{optequniav}
w+u=Tu,
\end{equation}
where $w$ is constant, and $u$ is an unknown function, which can be presented in multiple forms.

Policy iterations were introduced by Howard~\cite{Ho} for average-cost MDPs with finite state and action sets in two forms: for general problems and for unichain MDPs.  An MDP is unichain if every deterministic policy defines a Markov chain with one ergodic class.  For unichain MDPs the objective function $w^\phi$ is constant for all $\phi\in\F.$ Therefore, equation \eqref{optequniav} holds for unichain MDPs.

The numbers of variables and equations in LPs for unichain MDPs twice smaller than for general MDPs. However, Tsitsiklis~\cite{Tsi} proved that detecting whether an MDP is unichain is an NP-hard problem.  Thus, for a problem with unstructured data, it can be easier to find an average-cost optimal policy than to detect whether the problem is unichain or not.  For deterministic MDPs, detecting whether an MDP is unichain is a strongly polynomial problem~\cite{FYa}.
\section{Discounted MDPs with infinite sets of states and actions}\label{sec3}
Works by Blackwell~\cite{Blac1,Blac2} and Srauch~\cite{St} on discounted, positive, and negative MDPs with Borel states and actions were important contributions. Transition probabilities and costs are assumed to be Borel-measurable.   In addition, Blackwell~\cite{Blac1} provided an example showing that optimal values may not be Borel-measurable and $\epsilon$-optimal policies may not exist.  Blackwell, Freedman, Orkin  \cite{BFO,Freed} studied more general classes of models and policies. Bertsekas and Shreve~\cite{Bert} developed the theory for MDPs with Borel state and action sets and with universally measurable policies.  For countable-state MDPs, the theory of convergent MDPs was developed in \cite{FeSo,Fein1986,Fein1987}.  Convergent MDPs are more general than discounted, positive, and negative MDPs.

If certain continuity and compactness conditions hold for transition probabilities, cost functions, and action sets, then there exist deterministic optimal policies, which are defined by optimality conditions in the similar war as in the case of finite state and action sets, and the value functions $v_\beta$ can be computed by value iterations.  These conditions were formulated by Sch\"al \cite{Sch1975} for in the form of provided below Assumptions (S) and (W) for MDPs with setwise and weakly continuous transition probabilities respectively.

It is natural to consider MDPs with state spaces $\X$ and $\A$ being Borel subsets of Polish spaces.  In addition, for each state $x\in\X$ there  is a nonempty set $A(x)\subset\A$ of feasible actions.  It is assumed that the set $Gr_\X(A):=\{(x,a):\, x\in\X,\ a\in A(x)\}$  of feasible station-action pairs is a Borel subset of $\X\times\A$, and there exists a Borel mapping $\phi:\X\to\A$ such that $\phi(x)\in A(x)$ for all $x\in\X.$  These mappings are deterministic policies, and the set of deterministic policies is denoted by $\F.$  An arbitrary policy $\pi$ satisfies the property $\pi(A(x_t))|x_0,a),\ldots,x_t)=1$ for all histories from the set $H_t,$ up to each $t=0,1,\ldots.$

{\bf Assumption (W)}. [Sch\"al~\cite{Sch1975,Sch1993}]

\emph{(i) The set-value mapping $A:\X\to 2^\A$ is compact-valued and upper semicontinuous; }

\emph{(ii) the function $c(x,a)$ is lower semicontinuous  and bounded below on $Gr_\X (A);$}

\emph{(iii) the transition probability $p(\cdot|x,a)$ is weakly continuous  on $Gr_\X(A);$ that is, if $(x^{(k)},a^{(k)})\to (x,a)\in Gr_\X(A)$ and $ (x^{(k)},a^{(k)})\in Gr_\X(A),$ then for any bounded continuous function $f:\X\to\mathbb{R}$
\[\int_\X f(z)p(dz|x^{(k)},a^{(k)})\to \int_\X f(z)p(dz|x,a)\quad {\rm as}\quad k\to\infty.\]
}

{\bf Assumption (S)}. [Sch\"al~\cite{Sch1975,Sch1993}]

\emph{(i) The set-value mapping $A:\X\to 2^\A$ is compact-valued; }

\emph{(ii) the function $c$ bounded below on $Gr_\X (A),$ and for each $x\in\X$ the function $c(x,\cdot):\X\to\mathbb{R}$ is lower semicontinuous;}

\emph{(iii) for each $x\in\X,$ the transition probability $p(\cdot|x,a)$ is setwise continuous in $a$ on $A(x);$ that is, for each $x\in \X,$ if $a^{(k)}\to a\in A(x)$ and $ a^{(k)}\in A(x),$ then for any bounded continuous function $f:\X\to\mathbb{R}$
\[\int_\X f(z)p(dz|x,a^{(k)})\to \int_\X f(z)p(dz|x,a)\quad {\rm as}\quad k\to\infty.\]
}

As proved in Sch\"al~\cite{Sch1975,Sch1993}, each of Assumptions (W) and (S) implies the validity of optimality equations for discounted MDPs, these equations define optimal policies, and value iterations converge to optimal values.  In general, assumption (W) is more natural, and it is important for problems with incomplete information.  Assumption (S) does not require continuity of one-step costs and transition probabilities in the state variable.  Assumption (S) holds for MDPs with finite action sets and bounded below one-step costs without any additional assumption.

The proof under Assumption (W) is based on Berge's theorem which implies lower semicontinuity of the value function and upper semi-continuity of the solution multifunction for an optimization problem.  For topological spaces $U$ and $V,$ for a lower semicontinuous function $f:U\times V\to\mathbb{R},$ and for an upper semicontinuous set-valued function $F:\U\to 2^V$ with nonempty compact image set $F(u)$ for all $u\in U ,$ this theorem claims that  $f^*(u):=\min_{v\in F(u)} f(u,v)$ for all $u\in U,$ this function is lower semicontinuous, and the solution multifunction $F^*(u):=\{v\in F(v):\, f^*(u)=f(u,v)\}$ is compact-valued and upper semicontinuous.

In many models in operations research, including inventory control problems, action sets may not be compact, and the natural condition for one-step costs $c(x,a)$ is that for each $x\in\X$ the function $c(x,\cdot): A(x)\to\mathbb{R}$ is inf-compact, that is, for each $x\in \X$ and for each $\lambda\in\mathbb{R},$ the set $A_\lambda(x):=\{a\in A(x):\, c(x,a)\le \lambda\}$ is compact.

Luque-V\'asquez and  Hern\'andez-Lerma~\cite{LVHL} constructed an example showing that, if the assumption that the sets $F(u)$ are compact is replaced with the assumption that the lower semicontinuous function $f$ is inf-compact in the variable $v,$ then the conclusions of Berge's theorem fail.  This created the problem for extending the theory of MDPs to noncompact action sets.

In \cite{FKV,FKZ,FKZ1} this problem was resolved by introducing the class of $\K$-inf-compact functions for metric spaces $U$ and $V.$ A function $f:U\times V\to\mathbb{R}$  is called $\K$-inf-compact on $Gr_U(F)$, where $F:V\to 2^V\setminus\{\emptyset\},$ if for each compact $K\subset V$ the function $f:K\times V$ is inf-compact on $Gr_K(F).$

Many natural functions are $\K$-inf-compact.  For example, the function $f(u,v)=|u-v|$ and $f(u,v)=(u-v)^2$ are $\K$-inf-compact on $\mathbb{R}\times\mathbb{R}.$ Cost functions for inventory control problems are $\K$-inf-compact.  A function $c$ satisfying Assumptions (W)(i,ii) is $\K$-inf-compact on $Gr_\X(A).$  The following two assumptions generalize Assumptions (W) and (S) by expanding them to possibly noncompact action sets,

{\bf Assumption (${\rm W^*}$)}. \cite{FKZ}

\emph{(i) the function $c(x,a)$ is $\K$-inf-compact  and bounded below on $Gr_\X (A);$}

\emph{(ii) the transition probability $p(\cdot|x,a)$ is weakly continuous  on $Gr_\X(A).$
}

{\bf Assumption (${\rm S^*}$)}.\cite{FK}

\emph{(i) the function $c$ bounded below on $Gr_\X (A),$ and for each $x\in\X$ the function $c(x,\cdot):A(x)\to\mathbb{R}$ is inf-compact;}

\emph{(iii) for each $x\in\X,$ the transition probability $p(\cdot|x,a)$ is setwise continuous in $a$ on $A(x).$
}

Under Assumption~(${\rm S^*}$) the validity of the optimality equation, which defines an optimal  policy, and convergence of value iterations follow from a measurable selection theorem.   Hern\'andez-Lerma and Lasserre~\cite[Appendix D]{HLerma1} provided several formulations of measurable selection theorems useful for MDPs with compact action sets.  A selection theorem useful for noncompact action sets is introduced in \cite{FK}, and  is used there to prove that  Assumption~(${\rm S^*}$) is sufficient for the validity optimality equations, existence of deterministic optimal policies, and convergence of value iterations.

\section{Average-Cost MDPs with Infinite Sets of States and Actions}
\label{sec:4}
Average costs is a more difficult criterion to deal with than expected total discounted costs. There are many publications on this topic.  Arapostathis et al.~\cite{Ari} surveyed the results more than 30 years ago.  The survey has 208 citations and, according to Google Scholar,  was cited more than 700 times by October 2024. So, we are not trying to provide a comprehensive survey here.
\subsection{Finite-state Average-Cost MDPs with Infinite Action Sets}
We start with finite state MDPs. Let $X$ be finite.  The apparently natural assumptions for the existence of optimal policies is that the action sets $A(x)$ are compact, cost functions $c(x,a)$ and transition probabilities $p(z|x,a)$ are continuous in $a\in A(x)$ for all $x,z\in\X.$ This is true if the model is unichain \cite{Fein1975} or if every deterministic policy defines a Markov chain with the same number of recurrent classes. Let us consider  finite-dimensional sets of transition probabilities $P(x)=\{p(\cdot|x,a):\, a\in A(x)\}.$ These sets are compact if the sets $A(x)$ are compact and transition probabilities $p(z|x,a)$ are continuous in $a\in A(x).$ If each set $P(x)$ has a finite set of extreme points, then a deterministic optimal policy exist if all action sets $A(x)$ are compact, transition probabilities $p(z|x,a)$ are continuous in $a\in A(x),$ and one-step costs $c(x,a)$ are lower semicontinuous in $a\in A(x);$ see \cite[Theorem 2]{Fein1975}.
In general, Bather~\cite{Bath}, Chitashvili~\cite{Chit}, and Dynkin and Yushkevich~\cite{DY} provided examples demonstrating that optimal policies may not exist for MDPs with compact action sets and continuous one-step costs and transition probabilities.

Chitashvili~\cite{Chit} proved the existence of deterministic $\epsilon$-optimal policies. For the criterion
\[
w_1^\pi(x):=\liminf_{T\to\infty}\frac{1}{T} E_x^\pi\sum_{t=0}^{T-1} c(x_t,a_t)
\]
the existence of deterministic $\epsilon$-optimal policies. was proved in \cite{Fein1978}, and Bierth~\cite{Bi} extended these results to the criteria
\[
w_2^\pi(x):=E_x^\pi \limsup_{T\to\infty}\frac{1}{T} \sum_{t=0}^{T-1} c(x_t,a_t)
\]
and
\[
w_3^\pi(x):=E_x^\pi\liminf_{T\to\infty}\frac{1}{T} \sum_{t=0}^{T-1} c(x_t,a_t).
\]
 For communicating MDPs with compact action sets and continuous one-step costs, Bather~\cite{Bath1} proved the existence of deterministic optimal policies and the validity of the optimality equation.

For finite-state MDPs with arbitrary action sets and transition and cost functions, for every $\epsilon>0$ there exists an $\epsilon$-optimal Markov policy $\sigma$ such that $w^\sigma=w_1^\sigma$; \cite{Fein1980}.  This results implies that $w=w_1$ for MDPs with finite action sets.  Bierth~\cite{Bi} strengthened these results by showing that there exists an $\epsilon$-optimal Markov policy $\sigma$ such that   $w^\sigma=w_1^\sigma=w_2^\sigma=w_3^\sigma.$ This result implies that  $w=w_1=w_2=w_3$ if $\X$ is a finite set.  We recall here that $w(x)$ and $w_i(x),$ where $x\in\X$ and $i=1,2,3,$ are infimums of the corresponding objective criteria $w^\pi(x)$ and $w_i^\pi(x)$ over the set of all policies $\pi\in\Pi.$

\subsection{Average-Cost MDPs with Countable State Sets}
Canonical equations imply the existence of deterministic optimal policies under the assumption that the function $u$ is bounded \cite{DY,Pu,Ross68,Ross68a,Taylor}.  However, in many applications one-step costs costs $c$ are not bounded, and this typically implies that the function $u$ is unbounded.  Sennott~\cite{Senb,Sens} developed the theory for countable-state MDPs with unbounded costs.  These results can be viewed as generalizations of Bather's~\cite{Bath1} results to countable state spaces.  For such models the function $w$ is a constant, and only the second canonical equation should be considered.  The important observation was that that the second canonical equations can be replaced with an inequality.  Cavazos-Cadena~\cite{CC} provided an example of a countable-state MDP, for which the optimality inequality holds while the optimality equality does not hold.  Currently the major results  on countable-state MDPs follow from the available results on MDPs with Borel state spaces.

\subsection{Average-Cost MDPs with Borel State Spaces}
In this paper we follow the following terminology.  A Borel space is a Borel subset of a Polish space. A Polish space is a complete separable metric space. Let the state space $\X$ and action space $\A$ be Borel spaces.

Sch\"al \cite{Sch1993} considered the following assumption.

\textbf{Assumption (${\rm \bf G}$)}. $ w^*:=\inf\limits_{x\in
\X}w(x)<+\infty$.

If this assumption does not hold, then the problem is trivial because $w^\phi(x)=+\infty$ for every policy $\pi$ and for every state $x\in\X.$
Following Sch\"al~\cite{Sch1993}, define
the following quantities for a discount factor $\alpha\in [0,1)$:
\begin{equation*}
m_{\alpha}=\inf\limits_{x\in \X}v_{\alpha}(x),\quad
u_{\alpha}(x)=v_{\alpha}(x)-m_{\alpha},
\end{equation*}
\begin{equation*}
\underline{w}=\ilim\limits_{\alpha\uparrow
1}(1-\alpha)m_{\alpha},\quad\overline{w}=\slim\limits_{\alpha\uparrow
1}(1-\alpha)m_{\alpha}.
\end{equation*}
Observe that $u_\alpha(x)\ge 0$ for all $x\in\X.$ According to
Sch\"al \cite[Lemma 1.2]{Sch1993}, Assumption~(${\rm \bf G}$) implies
\begin{equation}\label{eq:schal} 0\le \underline{w}\le \overline{w}\le w^*<
+\infty.
\end{equation}

According to Sch\"al~\cite[Proposition 1.3]{Sch1993}, under
Assumption~(${\rm \bf G}$),  if there exists a Borel measurable function
$u:\X\to [0,+\infty)$ and a deterministic policy $\phi$ such that
\begin{equation}\label{eq7}
\underline{w}+u(x)\ge c(x, \phi(x))+\int_\X  u(y)q(dy|x,
\phi(x)),\quad x\in \X,
\end{equation}
then $\phi$ is \textit{average-cost optimal} and
$w(x)=w^*=\underline{w}=\overline{w}$ for all $x\in \X.$

Let us consider the following assumption introduced in \cite{Sch1993}.

\textbf{Assumption (B).} \emph{(i) Assumption~(${\rm \bf G}$)
holds, and (ii) $\sup_{\alpha\in [0,1)}u_{\alpha}(x)<\infty$ for all
$x\in \X$.}

 As proved by Sch\"al~\cite{Sch1993},  inequality \eqref{eq7} holds if Assumption (B)
 is added either to Assumption~(W) or to Assumption~(S).  In \cite{FKZ} it was shown that Assumption~(W)  can be replaced with more general Assumption~(${\rm W^*}$), which does not require compactness of action sets.  In \cite{FK} it was shown that  Assumption~(S)   can be replaced with the more general Assumption (${\rm S^*}$), which also does not assume compactness of action sets.  This fact was formulated in \cite{HL}, but the proof in \cite{HL} required a selection theorem proved in \cite{FK}.

Formula~\eqref{eq7} is called an optimality inequality. Another optimality inequality was introduced in \cite[Theorem 1]{FKZ}.   If Assumption~(${\rm \bf G}$) holds and there exists a Borel measurable
function $u:\X\to [0,+\infty)$ and a deterministic policy $\phi$ such
that
\begin{equation}\label{eq7111}
\overline{w}+u(x)\ge c(x, \phi(x))+\int_\X  u(y)q(dy|x,
\phi(x)),\quad x\in \X,
\end{equation}
then $\phi$ is average-cost optimal and
\begin{equation}\label{eq:7121}
w(x)=w^{\phi}(x)=\slim\limits_{\alpha\uparrow
1}(1-\alpha)v_{\alpha}(x)=\overline{w}=w^*,\quad x\in \X.
\end{equation}

The following weaker form of Assumption (B) is introduced in \cite{FKZ}.

\textbf{Assumption (${\rm  \underline{B}}$).}\emph{ (i)
Assumption~(${\rm \bf G}$)
 holds, and (ii) $\ilim\limits_{\alpha \uparrow
1}u_{\alpha}(x)<\infty$ for all $x\in \X$.}

\cite[Theorem 4]{FKZ}  states that optimality inequality \eqref {eq7111} holds for some measurable nonnegative measurable function $u$ if Assumptions~(${\rm W^*}$) and (${\rm \ \underline{B}}$) hold. \cite[Theorem 3.3]{FK} states the same conclusions under assumptions Assumptions~(${\rm S^*}$) and (${\rm  \underline{B}}$). \cite[Example 4.1]{FK} provides a countable-state Markov chain with costs satisfying Assumption (${\rm  \underline{B}}$) and is not satisfying Assumption (B).

\section{Markov Decision Processes with Incomplete Observations} Shiryaev~\cite{Sh1,Sh2} introduced and discussed many important results and ideas for control of stochastic processes with complete and incomplete observations in discrete and continuous time.  One of them is that the control problem can be reduced to the problem with states being probability distributions of the states of the unobserved model.  In statistics, these distributions are called prior or posterior depending on a concrete situation.  In control theory they are called belief states.  The reduction of MDPs with incomplete observations to MDPs with belief states was also described by Aoki~\cite{Ao}, {\AA}str\"em~\cite{As}, and Dynkin~\cite{Dy}, and currently this is the main method for studying and solving MDPs with incomplete observations.

Currently, the most popular model of an MDP with complete observations is a Partially Observable Markov Decision Process (POMDP), which broadly speaking is defined by the tuple $(\X,\Y, \A,\T, Q,c),$ where $\X$ is the space of hidden states, $\Y$ is the space of observations, $\A$ is the space of controls, $\T$ is the transition probability of the hidden process, $Q$ is the observation kernel, and $c$ is the one-step cost function.  Here $\X,$ $\Y,$ and $\A$ are Borel subsets of Polish spaces, $\T$ is a transition probability from $\X\times\A$ to $\X,$  $Q$ is the transition probability from $\A\times\X$ to $\Y,$  and $c:\X\times\A)\to \mathbb{R}\cup\{+\infty\}$ is a bounded below Borel function.

A POMDP can be reduced to a belief MDP $(\P(\X), \A, \bar{p}, \bar{c}),$ where $\P(\X)$ is the set of probability measures on the state space, $\bar{p}$ and $\bar{c}$ are the properly defined transition probability from $\P(\X)\times\A$ to $\P(\X)$ and one-step costs for the belief MDP; see \cite{HL} or \cite{FKZg,FKZg1} for details.  Here we use the notation $\P(E)$ for the set of probability measures on a measurable space $E.$ If $E$ is a Polish space, then $\P(E)$ is endowed with the topology of weak convergence of probability measures, and $\P(E)$ is also a Polish space.

Earlier studies, e.g.,  Rhenius~\cite{Rh} and  Yushkevich~\cite{Yu},  considered a more general model of a Markov Decision Process with Incomplete Information (MDPII) defined by a tuple $(\X,\Y, \A, P,c),$ where $\X,$ $\Y,$ and $\A,$ have the same meanings as for a POMDP, and $P$ is the transition probability from $\X\times\Y\times\A$ to  $\X\times\Y,$ and, for an MDPII, $c:\X\times \Y\times\A\to\R$ is the one-step cost function.  An MDPII can be reduced to the belief MDP $(\P(\X)\times\Y,\A,q,\bar{c}),$ where $q$ is the transition probability from $\P(\X)\times\Y\times\A$ to $\P(\X)\times\Y$ and $\bar{c}:\P(\X)\times\Y\times\A\to\R$ is the one-step cost function; see \cite{FKZg1} for details, where a POMDP is denoted as a ${\rm POMDP}_2.$

 The relation between  MDPIIs and POMDPs is that a POMDP is an MDPII with a specially defined transition probability $P.$  If the transitions for the MDPII are defined by the probability $P(dx_{t+1},dy_{t+1}|x_t,y_t,a_t),$ where $x_t,$ $y_t,$ and $a_t$ are the hidden state, observation, and control respectively at the epoch $t,$  then for the POMDP this probability is $P(dx_{t+1},dy_{t+1}|x_t,y_t,a_t)=Q(dy_{t+1}|a_t,x_{t+1})\T(x_{t+1}|x_t,a_t).$

 Rhenius~\cite{Rh} and Yushkevich~\cite{Yu} proved  the reduction of MDPIIs to belief MDPs for problems with Borel spaces $\X,$ $\Y,$ and $\A$ of hidden states, observations, and controls respectively. This result also implies the reduction of POMDPs to belief MDPs with the state space $\P(\X).$   An important question is whether  optimal policies exist for belief MDPs, and, if optimal policies exist, how to find them.  For discounted problems, according to \cite{FKZ}, if  Assumption~(${\rm W^*}$) holds  for the belief MDP, that is, for  the belief MDP the transition probability is weakly continuous and one-step cost is $\K$-inf-compact, then there are deterministic optimal policies for the belief MDP, and they can be computed by value iterations.  $K$-inf-compactness of the one-step cost function $\bar{c}$ for the belief  MDP follows from $\K$-inf-compactness of the original cost function $c;$ \cite[Theorem 3.3]{FKZg}. This is also true for MDPIIs. However, it is more difficult to verify weak continuity of the transition probability for the belief MDP.  For example, weak continuity of transition and observation probabilities is not sufficient for weak continuity of the transition probability for the belief MDP; \cite[Example 4.1]{FKZg}.

 For POMDPs, some conditions for weak continuity of transition probabilities for belief MDPs are given in monographs by Hern\'andez-Lerma~\cite{HL} and by Runggaldier and Stettner~\cite{RS}. According to \cite{FKZ},  weak continuity of the transition probability $\T$ and continuity of the observation probability $Q$ in total variation imply weak continuity of the transition probability $\bar{p}$ for the belief MDP. Another proof of this fact is given in  Kara, Saldi, Yuksel~\cite{KSY}, where it is also proved that continuity of the transition probability $\T$ in total variation implies weak continuity of the transition probability $\bar{p}$ for the belief MDP if the observation probability $Q$ does not depend on the control parameter $a.$  However, if $Q$ depends on $a,$ continuity of the transition probability $\T$ in total variation and continuity in total variation of the observation probability $Q$ in parameter $a$ imply weak continuity of the transition kernel $\bar{p}$ for the belief MDP \cite{FKZg}.  These results are summarized in the following theorem:

 \begin{theorem}[{\cite[Corollary 6.11]{FKZg1}}] \label{thm:feller:semi-uniform_feller}
   For a POMDP with the transition probability $\T$ and observation probability $Q,$  each of the following two conditions is sufficient for weak continuity of the transition probability $\bar{p}$ for the belief MDP:
    \begin{enumerate}
        \item[\rm{(i)}] $\T$ is weakly continuous, and $Q$ is continuous in total variation; \label{item:feller:transitions_weak_observations_total_variation}
        \item[\rm{(ii)}] $\T$ is continuous in total variation, and $Q$ is continuous in $a$ in total variation. \label{item:feller:transitions_total_variation_observations_total_variation_a}
    \end{enumerate}
\end{theorem}

%

 Continuity of transition probabilities for belief MDPs corresponding to MDPIIs was studied in \cite{FKZg,FKZg1}, where  the following definition was introduced. For Borel subsets $\S_1,$ $\S_2,$ and $\S_3$  of metric spaces, and a transition probability  $\Psi$ from  $\S_1\times\S_2$  to $\S_3$ is called {\it semi-uniform Feller} if, for each sequence $\{s_3^{\left(n\right)}\}_{n=1,2,\ldots}\subset\S_3$ that converges to $s_3 \in \S_3$ and for each bounded continuous function $f$ on $\S_1,$
\begin{equation}\label{eq:equivWTV3}
\lim_{n\to\infty} \sup_{B\in \B(\S_2)} \left| \int_{\S_1} f(s_1) \Psi(ds_1,B|s_3^{\left(n\right)})-\int_{\S_1} f(s_1) \Psi(ds_1,B|s_3)\right|= 0.
\end{equation}
It is proved in \cite{FKZg} that the transition probability $P$ from $\X\times\Y\times\A$ to  $\X\times\Y$ is semi-uniform Feller if and only if the transition probability $q$ from $\P(\X)\times\Y\times\A$ to $\P(\X)\times\Y$   for the belief MDP is semi-uniform Feller.  Semi-uniform Feller transition probabilities are weakly continuous.  Thus, semi-uniform continuity of the transition probability $P$ for the MDPII implies weak continuity of the transition probability for the corresponding belief MDP.  In particular, for POMDPs this result implies all the results on weak continuity of transition probabilities for belief MDPs stated in the previous paragraph.

\section{Discrete-Time Stochastic Filtering}\label{sec6}
POMDPs model discrete-time stochastic filtering problems defined by stochastic equations
\begin{subequations}\label{eq:filters:model}
    \begin{align}
        x_{t+1} &= F(x_t, a_t, \xi_t), && x_t \in \Xb, \quad a_t \in \Ab, \quad \xi_t \in \Xc, \label{eq:filters:transitions} \\
        y_{t+1} &= G(a_t, x_{t+1}, \eta_{t+1}), && a_t \in \Ab, \quad x_{t+1} \in \Xb, 
        \quad \eta_{t+1} \in \Hc, \label{eq:filters:observations}
    \end{align}
\end{subequations}
where $\Xc$ and $\Hc$ are Borel subsets of Polish spaces, $F$ and $G$ are Borel measurable functions, and $(\xi_t)_{t=0}^\infty$ and $(\eta_t)_{t=0}^\infty$ are two independent sequences of iid random variables with distributions $\mu$ and $\eta$ respectively; see \cite{FIKK} for details.
 The transition probability $\T$ and observation probability $Q$ for the POMDP for stochastic sequences $x_t$ and $y_t$ in \eqref{eq:filters:model} are
 \begin{subequations} \label{eq:filters:model_kernels}
   \begin{align}
        \T(B|x, a) &= \int_{\Xc} \mathbf{1}\{F(x, a, \xi) \in B\} \ \mu(d\xi), && B \in \B (\Xb), \quad x \in \Xb, \quad a \in \Ab, \label{eq:filters:transition_kernel} \\
        Q(C|a, x) & =\int_{\Hc}\mathbf{1}\{G(a, x, \eta) \in C\} \ \nu(d\eta), && C \in \B (\Yb), \quad a \in \Ab, \quad x \in \Xb. \label{eq:filters:observation_kernel}
    \end{align}
\end{subequations}
If the goal is to minimize expected total discounted costs \eqref{eqdicd}, and the bounded below cost function $c:\X\times\A\to\mathbb{R}\cup\{+\infty\}$ is $\K$-inf-compact, then an optimal policy exists, and it can be found by value iterations applied to the belief MDP if the transition probability $\bar{p}$ for the belief MDP is weakly continuous.  As explained above, this continuity depends on weak continuity and continuity in total variation of the transition kernels $\T$ and $Q.$

Both formulae   \eqref{eq:filters:model_kernels} can be rewritten in the same generic form
\begin{align} \label{eq:notation:stochastic_kernel}
    \kappa(B|s_2) := \int_{\Omega} \mathbf{1}\{\phi(s_2,\omega) \in B\} \ p(d\omega), &&
    B \in \B(\S_1), \quad
    s_2 \in \S_2,
\end{align}
where $\Omega,$ $\S_1,$ and $\S_2$ are Borel subsets of Polish spaces, $\phi : \S_2 \times \Omega \to \S_1,$ is a Borel measurable function, and $p$ is a probability measure on $\Omega.$

The questions are under which conditions the transition probability $\kappa$ from $\S_2$ to $\S_1$ is weak continuous and under which conditions it is continuous in total variation.  To answer these questions let us recall the following classic definitions.
\begin{definition} \text{(Continuity in distribution, total variation, and probability)}\label{def:notation:phi_continuity}
   \emph{ Let  $\S_1,$ $\S_2,$ and $\Omega$ be Borel spaces, and let  $p$ be a probability measure on $(\Omega, \B(\Omega)).$ A Borel function $\phi : \S_2 \times \Omega \to \S_1$ is continuous
    \begin{enumerate}
        \item[{\rm(i)}] \label{item:notation:continuity_in_distribution} in distribution $p$  if the function $s_2 \mapsto \int_\Omega f(\phi(s_2,\omega))\ p(d\omega)$ is continuous on $\S_2$ for every bounded continuous function $f : \S_1 \to \Rb;$
        \item[{\rm(ii)}] \label{item:notation:continuity_in_total_variation} in total variation with respect to (wrt) $p$ if for each $s_2 \in \S_2,$
        \begin{equation}
            \lim_{s_2' \to s_2}
            \sup_{B \in \B(\S_1)} \left|
                \int_{\Omega} \mathbf{1}\{\phi(s_2', \omega) \in B\} \
                -
                \mathbf{1}\{\phi(s_2, \omega) \in B\} \ p(d\omega)
            \right| = 0;
           \end{equation}
        \item[{\rm(iii)}] \label{item:notation:continuity_in_probability} in probability $p$ if $\phi(s_2',\:\cdot\:) \xrightarrow{p} \phi(s_2,\:\cdot\:)$ as $s_2' \to s_2$ for each $s_2\in \S_2,$ that is, for each $s_2 \in \S_2$ and each $\varepsilon > 0,$
        \begin{equation}
            \lim_{s'_2 \to s_2} p(\{
                \omega\in \Omega \,:\, \rho_{\S_1}(\phi(s_2',\omega),\phi(s_2,\omega)) \ge \varepsilon
            \}) = 0.
        \end{equation}
    \end{enumerate}}
\end{definition}

Continuity in total variation is stronger than continuity in distribution.  Continuity in probability is also stronger than continuity in distribution.  The first obvious observation is that $\kappa$ is weakly continuous if and only if $\phi$ is continuous in distribution.  The assumption, that the function $\phi$ is continuous, is used in several papers on control, and this assumption is much stronger than the assumptions that $\phi$ is continuous in probability or in distribution.  The second obvious observation is that $\kappa$ is continuous in total variation if and only if $\phi$ is continuous in total variation.  However, the condition that $\phi$ is continuous in total variation looks abstract, and we would like to have particular sufficient conditions for continuity of a transition kernel in total variation, which are useful for stochastic filtering. To do this, we recall the following theorem proved by Aumann\cite{Au}.

\begin{theorem} \text{\rm(\cite[Lemma F]{Au},~\cite[Lemma 1.2]{GS})} \label{thm:kernels:aumann}
    Let $\S_1$ and $\S_2$ be Borel spaces, and let $\kappa$ be a stochastic kernel on $\S_1$ given $\S_2.$  Then there exists a Borel measurable function  $\phi : \S_2 \times [0,1] \to \S_1$ such that
    \begin{align} \label{eq:kernels:kernel_representation}
        \kappa(B|s_2) &= \int_0^1 \mathbf{1}\{\phi(s_2,\omega) \in B\} \ d\omega,
        &&
        B \in \B(\S_1).
    \end{align}
\end{theorem}
This theorem implies the following corollary.
\begin{corollary}\text{\rm(\cite[Corollart 5.2]{FIKK})} \label{cor:kernels:aumann_multidim}
    Let $\S_1$ and $\S_2$ be Borel spaces, and let $\kappa$ be a stochastic kernel on $\S_1$ given $\S_2.$ Then for each natural number $n$ there exists a Borel measurable function $\phi : \S_2 \times [0,1]^n \to  \S_1,$ where the Borel $\sigma$-algebra is considered on the unit box $[0,1]^n,$ such that 
    \begin{align} \label{eq:kernels:multidim}
        \kappa(B|s_2) &= \int_{[0,1]^n} \mathbf{1}\{\phi(s_2, \omega) \in B\} \ d\omega,
        && B \in \B(\S_1).
    \end{align} 
\end{corollary}
The main differences between \eqref{eq:notation:stochastic_kernel} and \eqref{eq:kernels:multidim} are that $\Omega$ is a given Borel space, and $p$ is a given probability measure on $\Omega$ in \eqref{eq:notation:stochastic_kernel},  while  $\Omega=[0,1]^n,$ and  $p$ is the Lebesgue measure on $[0,1]^n$ in \eqref{eq:kernels:multidim}. 
In many filtering applications, the states and observations belong to Euclidean spaces, and $\S_1=\R^n,$ where $n=1,2,\ldots,$ represents state and observation spaces $\X$ and $\Y$ in the definitions of $\T$ and $Q$ in \eqref{eq:filters:model_kernels}.

 Let $D_x g = \frac{\partial g}{\partial x}$ denote the Jacobian of  a differentiable function $g : \Rb^n \to \Rb^n.$ The following condition is considered in \cite{FIKK}.

 {\bf Diffeomorphic Condition}
   \emph{ For the metric space $\S_2,$ open set $\Omega \subset \Rb^n,$  and a function $\phi : \S_2 \times \Rb^n  \to \Rb^n,$ the following statements hold:
    \begin{enumerate}
        \item[{\rm(i)}] $\phi$ is continuous on $\S_2 \times \Omega;$
        \item[{\rm(ii)}] $D_\omega\phi(s_2,\omega)$ exists for all $s_2 \in \S_2$ and $\omega \in \Omega;$
        \item[{\rm(iii)}] the matrix $D_\omega \phi(s_2,\omega)$ is nonsingular for all $s_2 \in \S_2$ and $\omega \in \Omega;$
        \item[{\rm(iv)}] the function $(s_2, \omega) \mapsto D_\omega \phi(s_2,\omega)$ is continuous on $\S_2 \times \Omega;$
        \item[{\rm(v)}] for each $s_2 \in \S_2$ the function $\omega \mapsto \phi(s_2,\omega)$ is a one-to-one mapping of $\Omega$  onto $\phi(s_2,\Omega).$
    \end{enumerate}}
The following theorem provides particular conditions for continuity of a transition probability in total variation, which are useful for stochastic filtering. In Theorem~\ref{thm:kernels:total_variation} and in the rest of this paper we follow the following remark.  A measurable subset $B$ of a measurable space $S$ can also be considered as a measurable space.  If a measure $m$ is defined on $B,$ we always consider and denote by the same letter the extension of this measure on $S$ by setting $m(S\setminus B)=0.$

\begin{theorem} \text{\rm (\cite[Theorem 5.4(b)]{FIKK}).} \label{thm:kernels:total_variation}
Let $\Omega$  be an open subset of $\Rb^n,$ where $n$ is a fixed natural number;
$p$ be a probability measure on $\Omega$ such that $p \ll \lambda^{[n]},$   where $\lambda^{[n]}$ is the Lebesgue measure on $\Rb^n;$
 $\S_1=\Rb^n;$ $\S_2$ be a Borel subset of a Polish space; and
a function $\phi : \S_2 \times \Rb^n  \to \Rb^,$ satisfy the Diffeomorphic Condition.
Then the transition probability $\kappa$ defined in \eqref{eq:notation:stochastic_kernel} is continuous in total variation.
\end{theorem}

Theorems~\ref{thm:feller:semi-uniform_feller} and \ref{thm:kernels:total_variation}  imply that each of the following two conditions  (i) and (ii) is sufficient for weak continuity of the transition probability for the belief MDP for a problem defined by equations~\eqref{eq:filters:model}:
   \begin{enumerate}
        \item[\rm(i)] \label{item:applications:ssm_case_1} The following statements (a) and (b) hold:
        \begin{enumerate}
            \item[\rm(a)] \label{item:applications:ssm_case_1a} the function $((x,a), \xi) \mapsto F(x,a,\xi)$ is continuous in distribution $\mu;$
            \item[\rm(b)] \label{item:applications:ssm_case_1b} $\Yb = \Rb^m,$ where $m$ is a natural number, $\Hc$ is an open  subset of $\Rb^m,$ $\nu \ll \lambda^{[m]},$ and the function $((a,x),\eta) \mapsto G(a,x, \eta)$ satisfies the Diffeomorphic Condition;
        \end{enumerate}
        \item[\rm(ii)] The following statements (a) and (b) hold:
        \begin{enumerate}
        \item[\rm(a)]
        \label{item:applications:ssm_case_2} $\Xb= \Rb^d,$ $\Xc$ is an open  subset of $\Rb^d,$ $\mu \ll \lambda^{[d]},$ and the function $((x,a),\xi) \mapsto F(x,a,\xi)$ satisfies the Diffeomorphic Condition;
        \item[\rm(b)] either the function $G$ does not depend on the parameter $a,$ that is, $G(a,x,\eta)=G(x,\eta), $ or $\Yb = \Rb^m,$  $\Hc$ is an open  subset of $\Rb^m,$ $\nu \ll \lambda^{[m]},$ and for each fixed $x\in\X$ the function $(a,\eta) \mapsto G(a,x, \eta)$ satisfies the Diffeomorphic Condition.
        \end{enumerate}
    \end{enumerate}

   If the goal is to minimize the expexted total discounted costs \eqref{eqdicd}, and the one-step cost function is $\K$-inf-compact and bounded, each of conditions (i) and (ii) implies the existence of optimal policies, the validity of optimality equations, and convergence of value iterations for problem \eqref{eq:filters:model}.       Examples of applications of Theorem~\ref{thm:kernels:total_variation} to stochastic filtering are described Theorem~\ref{corapplic} and in \cite{FIKK}.

   Two of these applications deal with filtering problems with two popular noise models: additive and multiplicative.  For additive noise, the function $\phi$ considered in the Diffeomorphic condition is $\phi(s_2,\omega)= f(s_2)+\xi(\omega),$ where $f:\S_2\to\R^n$ and $\xi$ is an $n$-dimensional random variable, and Theorem~\ref{thm:kernels:total_variation} implies that the function $\phi$ is continuous in total variation if the function $f$ is continuous. For multiplicative noise, the function $\phi$  is $\phi(s_2,\omega)= \diag(\xi(\omega))f(s_2),$ where $\diag(r)$ is the diagonal matrix whose diagonal entries are formed by the vector $r,$ and Theorem~\ref{thm:kernels:total_variation} implies that the function $\phi$ is continuous in total variation if the function $f$ is continuous and $f_j(s_2)\ne 0$ for all $j=1,\ldots,n.$

   The following theorem  generalizes  \cite[Corollary 7.3]{FIKK} to the case when the observation functions $G$ may depend on controls $a_t,$ $t=0,1,\ldots.$ Though in many papers dealing with filtering the function $G$ depends only on states and  noises, in important applications dealing with tracking, the observation function $G$ also depends on chosen controls.

   \begin{theorem}\label{corapplic} For an POMDP with  transition and observation probabilities defined in \eqref{eq:filters:model_kernels}, each of the following conditions is sufficient for  weak continuity of the transition probability $\bar{p}$ for the belief MDP:

(a) {\rm (additive transition noise)}  $F(x_{t},a_t,\xi_{t}) = f(x_{t},a_t) + \xi_{t},$ where  $\Xb=\Xc=\Rb^d,$ $f : \Xb\times\Ab \to \Xb$ is a continuous function, and  $\mu \ll \lambda^{[d]}$, and the function $G$ is measurable, and for each $x\in\X$ the function $(a,\eta) \mapsto G(a,x, \eta)$ satisfies the Diffeomorphic Condition;

(b) {\rm (multiplicative transition noise)}   $F(x_{t},a_t,\xi_{t}) = \diag(\xi_{t}) f(x_{t},a_t),$ 
        where   $\Xb = \Xc = \Rb^d,$   $f : \Xb\times\Ab \to \Xb$ 
        is a continuous function such that $f_j(x_{t},a_t))\ne 0$ for all $(x_{t},a_t)\in\Xb\times\Ab$ and for all $j=1,\ldots,d,$ and  
         $\mu \ll \lambda^{[d]},$ and for each $x\in\X$ the function $(a,\eta) \mapsto G(a,x, \eta)$ satisfies the Diffeomorphic Condition;

(c) {\rm (additive observation noise)}   $\tilde{G}(a_t,x_{t+1},\eta_{t+1}) = g(a_t,x_{t+1}) + \eta_{t+1},$ where $\Yb = \Hc = \Rb^m,$       $g :\A\times \Xb \to \Yb$ is a continuous function, 
        $\nu \ll \lambda^{[m]},$  and the function $((x,a), \xi) \mapsto F(x,a,\xi)$ is continuous in distribution $\mu,$ which takes place, for example, if $F$ is continuous;

(d)  {\rm (multiplicative observation noise)}   $G(a_t,x_{t+1},\eta_{t+1}) = \diag(\eta_{t+1})g(a_t,x_{t+1}),$ 
    where $\Yb = \Hc = \Rb^m,$ $g :\A\times \Xb \to \Yb$ is a measurable function such that this function is continuous in variable $x_{t+1}$ and, in addition, $g_i(a_t, x_{t+1})\ne 0$ for all $(a_t,x_{t+1})\in \A\times \Xb$ and for all $i=1,\ldots,m,$  
    $\nu \ll \lambda^{[m]},$  and the function $((x,a), \xi) \mapsto F(x,a,\xi)$ is continuous in distribution $\mu.$
\end{theorem}
\begin{proof}
The proof of of Theorem~\ref{corapplic} is similar to the proof of \cite[Corollary 7.3]{FIKK}, and it follows from Theorems~\ref{thm:feller:semi-uniform_feller}, \ref{thm:kernels:total_variation} and from sufficient conditions (i,ii) for weak continuity of the transition probability $\bar{p}$ for the belief MDP stated above in this section.
\end{proof}

\section{Kolmogorov's Equations for Jump Markov Processes and their Application to Continuous-Time Jump Markov Processes}
  This section mentions some of Albert Shiryaev's recent results  relevant to continuous-time MDPs. Some of these results advanced the theory of Markov processes.

Komogorov~\cite{Kol} introduced backward and forward equations for continuous-time Markov chains with finite and countable state spaces and for diffusion processes.  Feller~\cite{Fel} studied Kolmogorov's equations for jump Markov processes with Borel state spaces and unbounded jump rates.  A few years later, Feller noticed problems with formulations of forward equations, published an addendum~\cite{Fel}, and these problems remained open for more than 60 years.  Feinberg, Mandava and Shiryaev~\cite{FMS} solved them and then developed in \cite{FMS1} conditions for the validity of Kolmogorov's equations for jump Markov processes under more general assumption on unboundedness of jump rates than the assumptions introduced by Feller~\cite{Fel}; see also \cite{FS,FS1}.

These results were applied by Feinberg, Mandava and Shiryaev~\cite{FMS2} to the theory of Continuous-Time Markov Decision Processes (CTMDP), which deals with optimization of jump stochastic processes.  Monographs~\cite{GHL,KR,PZ} are devoted in to this theory.   In \cite{FMS2} it was proved that under broad assumptions an arbitrary policy for a CTMDP can be replaced with a Markov policy with the same or better performance.

\section{Concluding Remarks}  This article describes deep impacts of three papers \cite{Sh1,Sh2,VS} published by Albert Shiryaev in 1962-67.  In particular, \cite{VS}  contains the proof of the existence of nonrandomized stationary optimal policies for MDPs with finite state and action sets, when the objective is to optimize average rewards or costs per unit time.  References \cite{Sh1, Sh2} introduced important approaches and results on optimal decisions for discrete and continuous time problems with incomplete information including the reduction of problems with incomplete state observations to problems with complete state observations and with states being posterior probability distributions of states of the original problems. These posterior distributions are sometimes called beliefs, and the corresponding MDPs with complete information are called belief MDPs.  Currently these foundational results and approaches are broadly used in stochastic control, reinforcement learning, and artificial intelligence.

This paper  surveys some contemporary results on MDPs with infinite state and action sets with complete and incomplete information and on continuous-time MDPs. In particular, there was a recent significant progress in solving the longstanding open problem on providing sufficient conditions for weak continuity of transition probabilities for belief MDPs for problems with incomplete information, and this progress led to discovering broad sufficient conditions for the existence of optimal policies and for convergence of value iteration algorithms for discrete-time nonlinear filtering problems.  These results described in section~\ref{sec6} indicate deep relations between two fields: reinforcement learning and stochastic filtering, and both fields can benefit from these relations.  Another important recent development is solving in \cite{FMS,FMS1} Feller's~\cite{Fel} problem on the structure of solutions of forward Kolmogorov's equations for jump Markov processes; see also \cite{FS,FS1}.  These results significantly advanced the theory of continuous-time jump Markov decision processes by showing in \cite{FMS2} that under broad conditions  a Markov policy with the same or better performance can be constructed for an arbitrary policy.

\textbf{
Acknowledgement. }
Some of the research reported in this  publication was partially supported  by the U.S. Office of Naval Research (ONR) under Grant N000142412608.

%


\input{references}

\end{document}

%% file: references.tex
%
%
%


%% file: FeinbergSAFinalA.bbl
\begin{thebibliography}{99.}%
%
%
%
%

\bibitem{Ao} Aoki, M. (1965) Optimal control of partially observable Markovian systems. \emph{J. Franklin Inst.} \emph{280} pp. 367--386.


\bibitem{Ari}  Arapostathis, A.,  Borkar, V.S., Fernandez-Gaucherand, E, Ghosh M.K.,
Marcus, S.I. (1993) Discrete time controlled Markov processes with average
cost criterion: a survey, \textit{SIAM J. Control Optim.}
{31}(2) 282--344.

\bibitem{As}
 {\AA}str\"om, K.J. (1965). Optimal control of Markov processes with incomplete state information. \emph{J. Math.
Anal. Appl.} {10} pp. 174--205.

\bibitem{Au}
Aumann, R.J. (1964) Mixed and behavior strategies in infinite exstensive games. \emph{Advances in Game Theory}
\emph{52} pp. 627--650.

\bibitem{Bath}
Bather, J. (1973) Optimal decision procedures for finite Markov
chains. {P}art {I}:
  Examples.
\textit{Adv. in Appl. Probab.} {5} 328--339.

\bibitem{Bath1}
Bather, J. (1973) Optimal decision procedures for finite Markov
chains. {P}art {II}:
  Communicating systems.
\textit{Adv. in Appl. Probab.} {5} 521--540.

\bibitem{Ber} Berge, E. (1963) \textit{Topological Spaces.}
Macmillan, New York.

\bibitem{Bert} Bertsekas, D.P.,  Shreve, S.E. (1996) \textit{Stochastic Optimal Control: The Discrete-Time
Case.} Athena Scientific, Belmont, MA.

\bibitem{BT}  Bertsekas, D.P.,  Tsitsiklis, J.N. (1996)  \emph{Neuro-Dynamic Programming,} Athena Scientific, Belmont, MA.

\bibitem{Bi}  Bierth, K.-J. (1987) An expected average rward criterion. \emph{Stochastic Processes and Applications} {26} pp. 133--140.

\bibitem{BFZ}  Bishop, C.J.,     Feinberg, E.A.,  Zhang J. (2014) Examples concerning Abel and Cesaro limits. \emph{J. Math. Anal.  Appl.} {420}, pp. 1654--1661.

\bibitem{Blac}
Blackwell, D. (1962) Discrete dynamic programming. \textit{Ann. Math.
Statist.} {33}(2) pp. 719--726.

\bibitem{Blac1}
Blackwell, D. (1965) Discounted dynamic programming. \textit{Ann. Math.
Statist.} {36} pp. 226--235.

\bibitem{Blac2} Blackwell, D. (1967)  Positive dynamic programming.  In \emph{Proceedings of the fifth Berkeley symposium on mathematical statististics and probability (Berkeley, CA, 21 June-18 July 1965), vol. I: Theory of statistics.} Edited by L.M. Le Cam and J. Neyman. University of California Press (Berkeley and Los Angeles), pp. 415--418.

\bibitem{BFO}  Blackwell, D., Freedman, D., and  Orkin, M. (1974) The optimal reward operator in dynamic
programming, \emph{Ann. Probability} {2} pp. 926--941.



\bibitem{CC} Cavazos-Cadena, R. (1991) A counterexample on the optimality
equation in Markov decision chains with the average cost criterion.
\textit{Systems \& Control Lett.} {16}(5) 387--392.


\bibitem{Chit}
Chitashvili, R.Y. (1975) A controlled finite {M}arkov chain with an
arbitrary set of
  decisions.
\textit{Theor. Probability Appl.} {20}(4) 839--847.


\bibitem{DF}
Denardo, E.V.,   Fox, B.L. (1968)
Multichain Markov renewal programs.
\emph{SIAM J. Appl. Math.} {15}(3) pp. 468-487.


\bibitem{Der}
Derman, C.  (1962)
 On sequential decisions and {M}arkov chains.
\textit{Management Sci}. {9}(1) 16--24.



\bibitem{Dy}  Dynkin, E.B. (1965) Controlled random sequences. \emph{Theory
Probab. Appl.} {10} pp. 1--14.

\bibitem{DY}
 Dynkin, E.B.,   Yushkevich, A.A. (1979) \textit{{C}ontrolled
{M}arkov {P}rocesses}. Springer-Verlag, New York.

\bibitem{FedS} Federgruen, A., Schweitzer, P.J. (1980) SuccessivA survey of asymptotic value iteration for undiscounted Markov decision problems.  In: R. Hartley, L.C. Thomas, D.J. White (eds.) \emph{Resent Development in Markov Decision Processes,} Academic Press, New York, NY, pp. 73--109.

\bibitem{Fein1975} Feinberg, E.A. (1975) On controlled finite state Markov processes with compact control sets, \emph{Theor. Probab. Appl.} {20}, pp. 856--862.

\bibitem{Fein1978} Feinberg, E.A. (1978) The existence of a stationary $\epsilon$-optimal policy for a finite Markov chain \emph{Theor. Probab. Appl.} {23}, pp. 297--313.

\bibitem{Fein1980}
Feinberg, E.A. (1980)  An $\epsilon$-optimal control of a finite
Markov chain. \textit{Theor. Probab. Appl.} {25}(1)
70--81.

\bibitem{Fein1986}   Feinberg, E.A. (1986) Sufficient classes of strategies in discrete dynamic programming. I: Decomposition of randomized strategies and imbedded models. \textit{Theor. Probab. Appl.} {31} pp. 478--493.

    \bibitem{Fein1987}  Feinberg, E.A. (1987) Sufficient classes of strategies in discrete dynamic programming. II: Locally stationary strategies. \textit{Theor. Probab. Appl.} {32} pp. 658--668.

\bibitem{FHe} Feinberg, E.A.,  He, G. (2020)  Complexity bounds for approximately solving discounted MDPs   by value iterations,  \textit{Operations
  Research Letters} 48(5): 545--548

\bibitem{FHu}  Feinberg, E.A., Huang, J. (2014) The value iteration algorithm is not strongly
polynomial for discounted dynamic programming, \emph{Oper. Res. Lett.} 42: 130--131.


\bibitem{FIKK} Feinberg, E.A., Ishizawa, S.,   Kasyanov, P.O., Kraemer, D.N. (2025) Continuity of filters for discrete-time control problems defined by explicit equations. 	\emph{SIAM J. Control Optim.} 63(3): 1709-1735.

\bibitem{FIKK1} Feinberg, E.A., Ishizawa, S.,   Kasyanov, P.O., Kraemer, D.N. (2024) Sufficient conditions for solving statistical filtering problems by dynamic programming.  \emph{Proceedings of 63rd IEEE Conference on Decision and Control, December 16-19, 2024, Milan, Italy},  pp. 4052--4057.


\bibitem{FK} Feinberg, E.A.,   Kasyanov, P.O. (2021) MDPs with setwise continuous transition probabilities. \emph{Oper. Res. Lett.} {49}(5), pp. 734--740.


\bibitem{FKV} Feinberg, E.A.,   Kasyanov P.O., and M. Voorneveld, M. (2013) Berge's maximum theorem for noncompact image
sets. J.  Math. Anal. Appl. {397}(1):255--259.


\bibitem{FKZ} Feinberg, E.A.,   Kasyanov, P.O.,
Zadoianchuk, N.V. (2012) Average-cost Markov decision processes with weakly continuous transition probabilities,
\emph{Math. Oper. Res.} {37}, pp. 591-607.

\bibitem{FKZ1} Feinberg, E.A.,   Kasyanov P.O.,
Zadoianchuk, N.V. (2013) Berge's theorem for noncompact image sets. \emph{J.  Math. Anal.  Appl.} {397}, pp. 255--259.

\bibitem{FKZg} Feinberg, E.A.,   Kasyanov P.O., Zgurovsky, M.Z. (2016) Partially observable total-cost Markov decision
processes with weakly continuous transition probabilities. \emph{Math. Oper. Res.} {41}(2):
656--681.

\bibitem{FKZg1} Feinberg, E.A.,   Kasyanov P.O., Zgurovsky, M.Z. (2022) Markov decision processes with incomplete information
and semi-uniform Feller transition probabilities. \emph{SIAM J. Control Optim.}
{60}(4):2488--2513.

\bibitem{FKZg2} Feinberg, E.A.,   Kasyanov P.O., Zgurovsky, M.Z. (2023) Semi-uniform Feller stochastic kernels. \emph{J. Theor. Probab.} {36}, pp. 2262--2283.

%

\bibitem{FYa}  Feinberg, E.A., Yang, F. (2008) On polynomial classification problems for Markov decision processes. \emph{Oper. Res. Lett.} {36} pp. 527--530.
q


 \bibitem{FMS} Feinberg, E.A., Mandava, M., Shiryaev, A.N. (2014) On solutions of Kolmogorov's equations for jump Markov processes. \emph{J.  Math. Anal.  Appl.} {411}, pp. 261--270.

 \bibitem{FMS1} Feinberg, E.A., Mandava, M., Shiryaev, A.N. (2022) Kolmogorov's equations for jump Markov processes with unbounded jump rates. \emph{Ann. Oper. Res.} {317}(2), pp. 587--604.

 \bibitem{FMS2} Feinberg, E.A., Mandava, M., Shiryaev, A.N. (2022) Sufficiency of Markov policies for continuous-time jump Markov decision processes. \emph{Math.  Oper. Res.} {47}(2), pp. 1266--1286.

 \bibitem{FS} Feinberg, E.A., Shiryaev, A.N. (2022) Kolmogorov's equations for jump Markov processes and their applications to control problems. \emph{Theory Probab.  Appl.} {66}(4), pp. 582--600, 2022

 \bibitem{FS1} Feinberg, E.A., Shiryaev, A.N. (2024) On forward and backward Kolmogorov equations for pure jump Markov processes and their generalizations,
 \textit{Theor. Probab. Appl. } {68}(4), pp. 643--656.


 \bibitem{FeSo}  Feinberg, E.A., Sonin, I.M.  (1983) Stationary and Markov policies in countable state dynamic programming. \emph{Lecture Notes in Math.} 1021 pp. 111--129.



 \bibitem{Fel} Feller, W. (1940) On the integro-differential equations of purely discontinuous Markoff processes,
\emph{Trans. Amer. Math. Soc.,} {48}, pp. 488--515; Errata, Trans. Amer. Math. Soc., 58 (1945), p. 474.

\bibitem{Freed}  Freedman, D. (1974) The optimal reward operator in special classes of dynamic programming
problems, \emph{Ann. Probability} {2}, pp. 942-949.

\bibitem{GS}
  Gikhman, I.I.,  Skorohod, A.V. (1979) \emph{Controlled Stochastic Processes.} Springer, New York, NY.


\bibitem{GHL}  Guo, X.,   Hern\'andez-Lerma, O. (2009) \emph{Continuous-Time Markov Decision Processes: Theory and Applications.} Springer-Verlag, Berlin, 2009.



\bibitem{HL}  Hern\'{a}ndez-Lerma, O. (1989) \textit{Adaptive Markov Control Processes,} Springer-Verlag, New York.



\bibitem{HLerma} Hern\'{a}ndez-Lerma, O. 1991. Averege
optimality in dynamic programming on Borel spaces - Unbounded costs
and controls. \textit{Systems \& Control Lett.} {17}(3) pp.
237--242.

\bibitem{HLerma1} Hern\'{a}ndez-Lerma, O.,  Lassere, J.B. (1996)
\textit{Discrete-Time Markov Control Processes: Basic Optimality
Criteria}. Springer, New York.

\bibitem{Ho} Howard, R.A. (1960) \emph{Dynamic Programming and Markov Processes.} John Wiley \& Sons, New York, NY.

\bibitem{Kal}  Kallenberg, L.C.M. (1983)  \emph{Linear Programming and Finite Markovian Control Problems.} Mathematical Centre Tract 148, Mathematical Centre, Amsterdam.

\bibitem{Kal1}  Kallenberg, L.C.M. (2002) Finite state and action MDPs. E.A. Feinberg, A. Shwartz, eds.
\textit{Handbook of Markov Decision Processes. Methods and
Applications.} Kluwer, Boston, pp. 21--87.

\bibitem{KSY} Kara, A.D., Saldi, N., Y\"uksel, S. (2019) Weak Feller property of non-linear filters. \emph{Systems $\&$ Control Letters} {134}, 104512.

\bibitem{KR}  Kitaev, M.Yu.,   Rykov, V.V. (1995) \emph{Controlled Queueing Systems.} CRC Press, Boca Raton.

\bibitem{KM}   Kitahara, T.,  Mizuno, S. (2014) A bound for the number of different basic solutions
generated by the simplex method, \emph{Math. Program.} {137} pp. 579--586.

\bibitem{Kol} Kolmogoroff A. (1931) \"Uber die analytischen Methoden in der Wahrscheinlichkeitsrechnung,
\emph{Math. Ann.}, 104 ), pp. 415--458; English transl.:
A.N. Kolmogorov, On analytical methods in probability theory, in Selected Works of
A.N. Kolmogorov, Vol. II: Probability Theory and Mathematical Statistics, Math. Appl.
(Soviet Ser.) 26, Kluwer Acad. Publ., Dordrecht, 1992, pp. 62--108.

\bibitem{LVHL}
Luque-V\'asquez, F.,  Hern\'andez-Lerma. O. (1995)  A counterexample
on the semicontinuity of minima. \textit{Proc. Amer. Math. Soc.}
 (10) 3175--3176.

\bibitem{PZ}  Piunovskiy, A.B.  Zhang, Y. (2020)\emph{ Continuous-Time Markov Decision Processes.} Springer Nature,
Switzerland.

\bibitem{PY}  Post, I.  Ye, Y. (2015) The simplex method is strongly polynomial for deterministic
Markov decision processes, \emph{Math. Oper. Res.} {40} pp. 859--868.

\bibitem{Pu} Puterman, M.L. (1994) \emph{Markov Decision Processes: Discrete Stochastic Dynamic Programming.} John Wiley \& Sons, New York, NY.

\bibitem{Rh}  Rhenius, D. (1974) Incomplete information in Markovian decision models. \emph{Ann. Statist.} {2}(6) pp. 1327--1334.


\bibitem{Ross68}
 Ross, S.~M. (1968) Non-discounted denumerable {M}arkovian decision
model. \textit{Ann. Math. Statist.} {39}(2) 412--424.


\bibitem{Ross68a}
Ross, S.~M.  (1968) Arbitrary state {M}arkovian decision processes.
\textit{Ann. Math. Statist.} {39}(6) 2118--2122.



\bibitem{RS} {Runggaldier, W.J.,   Stettner, L. (1994)
\emph{Approximations of Discrete Time Partially Observed Control Problems,}
Applied Mathematics Monographs CNR, Giardini Editori, Pisa.}

\bibitem{Sch1975}  Sch\"{a}l, M. (1975) Conditions for optimality and for the limit of $n$-stage optimal policies to be optimal. \emph{Z. Wahrscheinlichkeitstheor. Verw. Geb.} {32} pp. 179--196.

\bibitem{Sch1993} Sch\"{a}l, M. (1993) Average
optimality in dynamic programming with general state space.
\textit{Math. Oper. Res}. {18}(1) 163--172.

\bibitem{Scher} Scherrer, B. (2016) Improved and generalized upper bounds on the complexity of
policy iteration. \emph{Math. Oper. Res.} {41} pp. 758--774.


\bibitem{Senb}
Sennott, L.I.  (1999) \textit{Stochastic Dynamic Programming and the
Control of Queueing
  Systems}.
John Wiley and Sons,
  New York.

  \bibitem{Sens}
Sennott, L.I.   2002.   Average reward optimization theory for
denumerable state spaces. E.A. Feinberg, A. Shwartz, eds.
\textit{Handbook of Markov Decision Processes. Methods and
Applications.} Kluwer, Boston, pp. 153--172.


\bibitem{Shapley} Shapley, L.S. (1963) Stochastic games, \emph{Proc. Natl. Acad. USA}  {39} pp. 1095--1100.

%
\bibitem{Sh1}
Shiryaev, A.N.  On the theory of decision functions and control by an observation process with incomplete data, \textit{Transactions of the Third Prague
Conference on Information Theory, Statistical Decision Functions,
Random Processes} (Liblice, 1962), 1964, pp.~657-681 (in Russian); Engl.
transl. in \textit{Select. Transl. Math. Statist. Probab.}
6(1966), 162--188.

\bibitem{Sh2}
Shiryaev, A.N.  Some new results in the theory of controlled
random processes. \textit{Transactions of the Fourth Prague
Conference on Information Theory, Statistical Decision Functions,
Random Processes} (Prague, 1965), 1967, pp.~131-201 (in Russian); Engl.
transl. in \textit{Select. Transl. Math. Statist. Probab.}
8(1969), 49--130.
%



\bibitem{Silver}  Silver, D.,  Hubert, T.,   Schrittwieser, J.,  Antonoglou, I.,  Lai, M.,  Guez, A.,  Lanctot, M.,   Sifre, L.,   Kumaran, D.,  Graepel, T.,  Lillicrap, T.,  Simonyan, K.,  Hassabis D. (2018) A general reinforcement learning algorithm that masters chess, shogi, and Go through self-play, \emph{Science} 3662(6419) pp. 1140-1144.

\bibitem{St}    Strauch, R. (1966) Positive Dynamic Programming.  \emph{Ann. Math. Statist.} {37}  pp. 871--890.

\bibitem{SB} Sutton, R.S.,  Barto, A.G. (2018) \emph{Reinforcement Learning: An Introduction (2nd ed.).} The MIT Press, Cambridge, MA.

\bibitem{Taylor}
 Taylor, III, H.~M.. 1965. Markovian sequential replacement
processes. \textit{Ann. Math. Statist.} {36}(6) 1677--1694.

\bibitem{Ts}  Tseng, P. (1990)  Solving h-horizon, stationary Markov decision problems in time
proportional to log(h), \emph{Oper. Res. Lett.} {9} pp. 287--297.

\bibitem{Tsi} Tsitsiklis, J.N. (2007) \emph{NP}-hardness of checking the
unichain condition in average cost MDPs. \emph{Oper.
Res. Lett.} {35} pp. 319--323.


\bibitem{VS}
 Viskov, O.V.,  Shiryaev, A.N. On controls leading to optimal stationary regimes, Proceedings of the Steklov Institute of Mathematics, 71 (1964), pp. 35-45 (In Russian);  English translation: Report Number FTD-HT-67-69, National Technical Information Service, U.S. Department of Commerce.

\bibitem{Ye}  Ye, Y. (2011) The simplex and policy-iteration methods are strongly polynomial
for the Markov decision problem with a fixed discount rate, \emph{Math. Oper. Res.} \emph{36} pp. 593--603.

\bibitem{Yu} Yushkevich, A.A. (1976) Reduction of a controlled Markov model with incomplete data to a problem with complete information
in the case of Borel state and control spaces. \emph{Theory
Probab. Appl.} \emph{21}pp. 153--158.

\end{thebibliography}
